\documentclass[12pt,oneside]{article}
\usepackage{amsmath}
\usepackage{graphicx}
\usepackage{geometry}
\usepackage{amssymb, gensymb}
\usepackage{amsthm}
\usepackage{mathtools}
\usepackage{mathrsfs}
\usepackage{indentfirst}
\usepackage{color}
\usepackage{datetime}
\usepackage{bm}
\usepackage{geometry}
\usepackage[english]{babel}

\begin{document}

\begin{center}
    
{\bf Speed of convergence of Chernoff approximations\\ to solutions of evolution equations}

 A.\,V.~Vedenin, V.\,S.~Voevodkin, V.\,D.~Galkin, E.\,Yu.~Karatetskaya, I.\,D.~Remizov\footnote{ivremizov@yandex.ru}

National Research University Higher School of Economics

October 2019
\end{center}


Abstract. This communication is devoted to establishing the very first steps in study of the speed at which the error decreases while dealing with the based on the Chernoff theorem approximations to one-parameter semigroups that provide solutions to evolution equations.


\vskip5mm

\textbf{Introduction.} Since the middle of the XX century it is a well known fact \cite{remiz:1,remiz:2} that the solution of a well-posed Cauchy problem for a linear evolution partial differential equation (examples: Sch\"odinger equation, parabolic equations) is given by a strongly continuous semigroup of linear bounded operators whose infinitesimal generator is a (usually unbounded) linear operator from the right-hand side of the evolution equation. Let us explain this in more details and introduce some notation which will be useful for the main text.

Let $X$ be an infinite set, and $\mathcal{F}$ be a Banach space of (not necessarily all) number-valued functions on $X$, and let $L$ be a closed linear operator $L\colon Dom(L)\to\mathcal{F}$ with the domain $Dom(L)\subset\mathcal{F}$ dense in $\mathcal{F}$. We consider the Cauchy problem for the evolution equation
$$
\left\{ \begin{array}{ll}
u'_t(t,x)=Lu(t,x),\\
u(0,x)=u_0(x),
\end{array} \right.\eqno(1)
$$
where $x\in X$, $u_0\in\mathcal{F}$, $u(t,\cdot)\in \mathcal{F}$ for all $t\geq 0$, and $L$ is, for example, in trivial case the Laplace operator $\Delta$ (so $u'_t=Lu$ is the heat equation), or (in nontrivial case) a more sophisticated linear differential operator with variable coefficients that do not depend on $t$ but depend (usually nonlinearly) on $x$. It is known \cite{remiz:1,remiz:2} that in case of existence of the $C_0$-semigroup $\left(e^{tL}\right)_{t\geq 0}$ with the generator $(L, Dom(L))$ the solution to the Cauchy problem (1) exists (in sense that l.h.s. is equal to r.h.s. in $\mathcal{F}$) and is given by the equality $u(t,x)=(e^{tL}u_0)(x)$ for all $t\geq 0$ and $x\in X$. If  $u_0\in Dom(L)$, then $u(t,\cdot)\in Dom(L)$ for all $t\geq 0$ and the solution $u$ is a classical solution (in the terminology of \cite{remiz:1}), and for arbitrary $u_0\in\mathcal{F}$ the solution of the Cauchy problem exists only as the solution of the corresponding integral equation $u(t,\cdot)=L\int_0^tu(s,\cdot)ds + u_0$. We write sometimes $u(t,x)$ and sometimes  $u(t,\cdot)$ assuming that the role of $\mathcal{F}$ can be played by, for example, the $L^p(\mathbb{R})$ space, then (1) holds only for almost all $x\in\mathbb{R}$. We see that in this case the notation $u(t,x)$ is not completely precise because all the versions of the function $x\longmapsto u(t,x)$ correspond to the same vector $u(t,\cdot)\in L^p(\mathbb{R})$, which usually does not lead to misunderstanding.

Equality $u(t,x)=(e^{tL}u_0)(x)$ shows that finding the semigroup  $\left(e^{tL}\right)_{t\geq 0}$ is a hard problem because it is equivalent to solving the Cauchy problem (1) for each $u_0\in\mathcal{F}$. However, if the so-called Chernoff function is constructed, a function $G$ which satisfies the conditions of the Chernoff theorem (in particular, satisfies $G(t)=I+tL+o(t)$, $t\to+0$), then the semigroup is given by the equality $e^{tL}=\lim_{n\to\infty}G(t/n)^n$. An advantage of this approach arises from the fact that usually it is possible to define $G$ by an explicit and not very long formula which contains coefficients of operator $L$, thus obtaining approximations to the solution of the Cauchy problem (1) converging to the solution in $\mathcal{F}$ as $n\to\infty$. Expressions $G(t/n)^nu_0$ are called Chernoff approximations to the solution of the Cauchy problem (1).

This communication is dedicated to the study (for fixed $t>0$) of the speed of decreasing (depending on $f\in\mathcal{F}$) of the norm of the difference between approximate and exact solution $\|G(t/n)^nf - e^{tL}f\|$ as $n\to\infty$. First definitions in this new area of $C_0$-semigroups studies are given, model examples of semigroups and their Chernoff approximations are examined, and the first (model) example of fast-converging Chernoff approximations for approximate calculation of (known before the Chernoff theorem appearance) solutions of one-dimensional heat equation is presented.

\textbf{1. Approximation subspaces in the Chernoff theorem.} The definition of Chernoff tangency was introduced in \cite{remiz:4} and will play the key role in what follows. 

{\sc Definition 1.} We say that operator-valued function $G$ is \textit{Chernoff tangent} to the operator $L$ (the details are provided below) iff the following conditions (CT0)-(CT4) are fulfilled:

(CT0). Hereafter $\mathcal{F}$ is a Banach space, and $\mathscr{L}(\mathcal{F})$ is the space of all linear bounded operators in $\mathcal{F}$. Let $G$ be a mapping $G\colon [0, +\infty) \to \mathscr{L}(\mathcal{F})$, in other words a family of linear bounded operators $(G(t))_{t\geq 0}$.  A closed linear operator $L\colon Dom(L) \to \mathcal{F}$ has the domain $Dom(L)\subset\mathcal{F}$ dense in $\mathcal{F}$.

(CT1). The family $G$ is strongly continuous (= continuous in strong operator topology in the space $\mathscr{L}(\mathcal{F})$), i.~e. the mapping $t \longmapsto G(t)f\in\mathcal{F}$ is continuous on $[0, +\infty)$ for each $f \in \mathcal{F}$;

(CT2). $G(0) = I,$ i.~e. $G(0)f=f$ for each $f \in \mathcal{F}$;

(CT3). There exists such a dense in $\mathcal{F}$ linear subspace $\mathcal{D} \subset \mathcal{F}$, that for all $f \in \mathcal{D}$ there exists a limit $\lim_{t \to 0}(G(t)f-f)/t$. Let us denote the value of this limit as $G'(0)f$;

(CT4). The closure of the linear operator $(G'(0), \mathcal{D})$ exists and is equal to $(L, Dom(L))$.

{\sc Remark 1.} The fact that $Dom(L)$ is dense in $\mathcal{F}$ follows from (CT3) and (CT4), so actually we may omit that in (CT0). We now can state Classic Chernoff's theorem (Paul Robert Chernoff \cite{Chernoff}, 1968) in the following wording, with (E)xistence and (N)orm growth conditions being separate from (CT), cf. \cite{remiz:1,remiz:2, Chernoff,remiz:4}.

\textsc{The Chernoff theorem, contemporary wording.} (I.D.Remizov \cite{remiz:4}, 2016)
Let $\mathcal{F}$ be a Banach space. Suppose there is a mapping  $G\colon [0, +\infty) \to \mathscr{L}(\mathcal{F})$ and a closed linear operator $L\colon Dom(L) \to \mathcal{F}$ with the domain $Dom(L)\subset\mathcal{F}$. Suppose the following conditions hold:

(E). There exists a  $C_0$-semigroup  $(e^{tL})_{t \geq 0}$ with the infinitesimal generator $(L, Dom(L))$;

(CT). $G$ is Chernoff tangent to $L$;

(N). There exists such a number $\omega \in \mathbb{R}$ that $\|G(t)\| \leq e^{\omega t}$ for all $t\geq 0$.

Then for each $f \in \mathcal{F}$ and $T>0$ it is true that $$\lim_{n\to\infty}\sup_{t\in[0,T]}\left\|\left(G\left(\frac{t}{n}\right)\right)^nf - e^{tL}f\right\|=0.\eqno(2)$$

Suppose that due to the Chernoff theorem equality (2) is true, i.e. we have the convergence  $G(t/n)^nf\to e^{tL}f$ for each $f\in\mathcal{F}$. But what is the speed of that convergence, how fast the error decreases to zero as $n$ tends to infinity? Moreover, what questions on the convergence rate we can state, how to measure it, what to expect and what not to expect? Besides consideration of several particular cases \cite{OSS2012} the building of this theory is not finished yet, and more and more researchers are involved into this area \cite{Zag, Gom}. Let us present some preparatory ideas as a first step on that way. On the one hand, we can for each $t>0$ and each Chernoff function $G$ define a function $C_G(t)$ which maps the space $\mathcal{F}$ into the space  $c_0$ of sequences tending to zero by the rule
$(C_G(t)f)(n)=\|G(t/n)^nf-e^{tL}f\|$ and study its properties as they can give us full information on the subject that we are interested in. On the other hand, this function is nonlinear and has too many parameters/arguments ($G$, $t$, $f$) to study it easily. Nevertheless, everything that we will know on the convergence of Chernoff approximations will be a corresponding statement about this function.

Below we will use (as $n\to\infty$) standard symbols ''o small`` and ''O big``; recall that it follows from the definitions of that symbols that $a_n=o(b_n)$ implies $a_n=O(b_n)$.

\textsc{Proposition 1.} (I.D.Remizov \cite{Remizov-NNGU2018}, 2018) Suppose that $\tau \subset [0, +\infty)$, $\tau\setminus\{0\}\neq\varnothing$ and $w\colon [1,+\infty)\to [0,+\infty)$ such that $\lim\limits_{x\to+\infty}w(x)=0$. Then the set
$$A_w^\tau=
\Big\{f\in\mathcal{F} : \sup_{t\in\tau}\|G(t/n)^nf-e^{tL}f\|=O(w(n))\textrm{ as }n\to\infty\Big\}$$
is a linear subspace in $\mathcal{F}$. Moreover, it follows from $$w_2(x)=o(1), w_1(x)=O(w_2(x))\textrm{ as } x\to+\infty$$ that 
$A^\tau_{w_1}\subset A^\tau_{w_2}.$ (In other words, the error (when $n\to\infty$) decreases on vectors $f\in A^\tau_{w}$ not slower than  $\mathrm{const}\cdot w(n)$. Inclusion $A^\tau_{w_1}\subset A^\tau_{w_2}$ means that the error decreases on vectors from  $A^\tau_{w_1}$ not slower that on vectors from $A^\tau_{w_2}$.)

\textsc{Definition 2.} (I.D.Remizov \cite{Remizov-NNGU2018}, 2018) Let us call $A^\tau_w$ an approximation subspace of order $w(n)$ on the set $\tau$ for Chernoff function $G$, and let us call the inclusion $A^\tau_{w_1}\subset A^\tau_{w_2}$ a hierarchy of approximation subspaces associated with Chernoff function $G$. We call arbitrary linear subspace $K\subset\mathcal{F}$ an approximation subspace iff there exists such a function $w(n)\to 0$ that $A^\tau_w=K$. 

\textsc{Remark 3.} The proof of proposition 1 is a simple check: suppose that numbers $\alpha$ and $\beta$ are arbitrary and vectors $f$ and $g$ belong to $A^\tau_w$; let us prove that $h=\alpha f+\beta g$ also belongs to $A^\tau_w$:
$$\|G(t/n)^nh-e^{tL}h\|=
\|G(t/n)^n(\alpha f+\beta g)-e^{tL}(\alpha f+\beta g)\|$$
$$\leq|\alpha|\cdot\|G(t/n)^nf-e^{tL}f\|+|\beta|\cdot\|G(t/n)^ng-e^{tL}g\|=O(w(n))+O(w(n))=O(w(n)).$$
To finish the proof it suffices to take $\sup_{t\in\tau}$ from a l.h.s. and r.h.s. of the inequality. Inclusion $A^\tau_{w_1}\subset A^\tau_{w_2}$ follows directly from the definition of $A^\tau_w$. Also note that $\tau_1\subset\tau_2$ implies $A^{\tau_2}_{w}\subset A^{\tau_1}_{w}$ as the supremum of a non-negative function on a smaller set does not exceed the corresponding one taken on a bigger set. Let us now suppose that $\tau$ is fixed, so we will write $A_{w}$ instead of $A^\tau_{w}$.

\textsc{Remark 4.} The set $A_w$ is constructed uniquely having $w$ (this is obvious), but it is not possible to uniquely reconstruct $w$ having $A_w$: for example, there is the same approximation subspace corresponding to functions $n\longmapsto w(n)$ and $n\longmapsto 2w(n)+w(n)/n$ as each one is O big of another. 

\textsc{Remark 5.} Each vector $f\in\mathcal{F}$ belongs to some approximation subspace, in particular, to the space $A_w$ where $w(n)=\|G(t/n)^nf-e^{tL}f\|$. This is why we have $\mathcal{F}=\bigcup_{w(n)=o(1),w(n)\geq0}A_w$, i.e. each vector (with the approximation subspace that it belongs to) takes its place in a hierarchy, which possesses a natural structure of an ordered set (partially ordered in non-degenerate cases).

\textsc{Remark 6.} Not every linear subspace in $\mathcal{F}$ is an approximation subspace. For example, if we set $G(t)=e^{tL}$, then $\|G(t/n)^nf-e^{tL}f\|=0$ for all $f$, so for that Chernoff function $G$ there is only one approximation subspace, $\mathcal{F}$ itself. Moreover, $\mathcal{F}=A_w$ for each function $w(n)=o(1)$. This example also shows that the speed of convergence in the Chernoff theorem can be arbitrary high if the Chernoff function $G$ is chosen luckily.

\textsc{Remark 7.} In the paper \cite{Remizov-NNGU2018} the following conjecture was presented: the higher the number of common derivatives with respect to $t$ at zero functions $t\mapsto e^{tL}$ and $t\mapsto G(t)$ have, the higher the speed of convergence one can obtain on non-trivial subspaces. For $S(t)=I+tL+o(t)$ it is possible to have $C/n$ convergence speed, whereas for $S(t)=I+tL+\frac{1}{2}t^2L^2+o(t^2)$ it appears to have $C/n^2$ rate, and etc.

\textbf{2. Semigroup of translations on the real line.} Let us consider the case $X=\mathbb{R}$, $(Lf)(x)=f'(x)$, and let $\mathcal{F}=UC_b(\mathbb{R})$ be a Banach space of all bounded and uniformly continuous functions  $f\colon \mathbb{R}\to \mathbb{R}$ with the uniform norm $\|f\|=\sup_{x\in\mathbb{R}}|f(x)|$. Cauchy problem (1) in that case takes the form $[u'_t(t,x)=u'_x(t,x), u(0,x)=u_0(x)]$, and its solution is the function $u(t,x)=u_0(x+t)$, so $(e^{tL}f)(x)=f(x+t)$, i.e. the function $f$, ''translated`` on $t$. Direct checking shows that $(e^{tL})_{t\geq 0}$ fits the definition of a $C_0$-semigroup. This example appears to be rich enough to use it for answering some general questions.

\textsc{Theorem 1.} Speed of convergence in the Chernoff theorem can be arbitrary slow. That is, if the function $w$ is given, and  $\lim_{x\to+\infty}w(x)=0$, then there exists such $X, \mathcal{F}, L, e^{tL}, G, f$, that $w(n/t)=O\left(\|G(t/n)^nf-e^{tL}f\|\right)$ as $n\to\infty$ for each $t>0$.

\textsc{Proof.} Let us go back to the example of a translation semigroup and set $(G(t)f)(x)=f(x+t+tw(1/t))$ for $t>0$ and $G(0)f=f$. Checking the conditions (CT0)-(CT4) shows that $G$ is a Chernoff function for the semigroup of translations. It follows from the definition of function $G$ that  $(G(t/n)f)(x)=f(x+t/n+(t/n)w(n/t))$ and $(G(t/n)^nf)(x)=f(x+t+tw(n/t))$. If the function $f$ has a continuous derivative then  $(G(t/n)^nf)(x)-(e^{tL}f)(x)=f(x+t+tw(n/t))-f(x+t)=f'(\xi)tw(n/t)$ where $\xi\in[x+t,x+t+tw(n/t)]$. Hence $\|G(t/n)^nf-e^{tL}f\|=\sup_{x\in\mathbb{R}}|(G(t/n)^nf)(x)-(e^{tL}f)(x)|=\sup_{x\in\mathbb{R}}|f'(\xi)tw(n/t)|\geq |f'(\xi)|tw(n/t)$ 
for each $t,x,n$. Fix arbitrary $t>0$ and set $x=-t$, then it appears that $\xi\in[0,tw(n/t)]$, so $\lim_{n\to\infty}\xi=0$ as $\lim_{x\to+\infty}w(x)=0$. Therefore, for $f(x)=\sin(x)$ we have $\lim_{n\to\infty}|f'(\xi)|=\cos(0)=1$, and, in particular, there exists such a number $n_0$ that for all $n>n_0$ and $x=-t$ we have $|f'(\xi)|>\frac{1}{2}$. 
So for all $n>n_0$ we have $\|G(t/n)^nf-e^{tL}f\|\geq \frac{1}{2}tw(n/t)$, that is equivalent to $w(n/t)\leq \frac{2}{t}\|G(t/n)^nf-e^{tL}f\|$, which means that $w(n/t)=O\left(\|G(t/n)^nf-e^{tL}f\|\right)$ for $f(x)=\sin(x)$.  $\Box$

\textsc{Remark 8.} Let us consider the same semigroup of translations and Chernoff function $(G(t)f)(x)=f(x+t+at^2)$, where $a\neq0$. In that case $(G(t/n)^nf)(x)=f(x+t+at^2/n)$, and if $f$ is a H\"older function (i.e. $|f(x_1)-f(x_2)|\leq C |x_1-x_2|^\alpha$ for some $C\geq0$ and $0<\alpha\leq1$), then $\|G(t/n)^nf-e^{tL}f\|= \sup_{x\in\mathbb{R}}|(G(t/n)^nf)(x)-(e^{tL}f)(x)|=\sup_{x\in\mathbb{R}}|f(x+t+at^2/n)-f(x+t)|\leq C|at^2/n|^\alpha=Ct^{2\alpha}|a|^\alpha\left(\frac{1}{n}\right)^\alpha=O\left(\left (\frac{1}{n}\right)^\alpha\right)$, where we set $x_1=x+t+at^2/n$ and $x_2=x+t$. So each intersection of $UC_b(\mathbb{R})$ with the space of H\"older functions with coefficient $0<\alpha\leq1$ is and approximation subspace of order $w_\alpha(n)=\left(\frac{1}{n}\right)^\alpha$. Here we have found a fragment of approximation structure for Chernoff function $(G(t)f)(x)=f(x+t+at^2)$. It seems interesting to find this hierarchy completely and to examine the case when $a$ is dependent on $x$.

\textbf{3. Semigroup of solutions of heat equation on the real line.} Let it be again $X=\mathbb{R}$ and $\mathcal{F}=UC_b(\mathbb{R})$, but $(Lf)(x)=a^2f''(x)$ for fixed $a>0$. Cauchy problem (1) in that case takes form  $[u'_t(t,x)=a^2u''_{xx}(t,x), u(0,x)=u_0(x)]$ which is a heat equation whose solution (and hence the semigroup $(e^{tL})_{t\geq 0}$) are given by Poisson integral $u(t,x)=(e^{tL}u_0)(x)=\int_{\mathbb{R}}\Phi(x-y,t)u_0(y)dy$, where $\Phi(x,t)=(2a\sqrt{\pi t})^{-1}\exp\left(\frac{-x^2}{4a^2t}\right)$. Chernoff function for a more general equation with variable coefficients was found in \cite{Rem-add-2}, and in this particular case it takes form $(G(t)f)(x)=\frac{1}{4}f(x+2a\sqrt{t})+\frac{1}{4}f(x-2a\sqrt{t})+\frac{1}{2}f(x)$. In A.V.Vedenin's diploma thesis (Higher School of Economics, Nizhny Novgorod, 2019) another Chernoff function was constructed: $(S(t)f)(x)=\frac{2}{3} f\left(x \right)+\frac{1}{6} f\left(x+a\sqrt{6t} \right)+\frac{1}{6} f\left(x-a\sqrt{6t} \right)$. 

\textsc{Conjecture 1.} Chernoff function $S$ has an approximation subspace of order $1/n^2$, and for some vectors $u_0\in UC_b(\mathbb{R})$ it provides faster convergence to the exact solution of the heat equation than function $G$ does.

\textsc{Remark 9.} The discussion of the above-mentioned model examples rises a hope for the creation of universal methods of construction of fast-converging Chernoff approximations (in particular, Feynman formulas \cite{remiz:3} and their analogues \cite{R-Starod, Loboda}) for evolution equations with variable coefficients. 

\textbf{Acknowledgments.}  The research is supported by Laboratory of Dynamical Systems and Applications NRU HSE, of the Ministry of science and higher education of the RF grant ag. No 075-15-2019-1931. I.D.Remizov is also grateful to Prof. O.G.Smolyanov for moral support during this research, to D.V.Turaev for fruitful discussions and to O.E.Galkin and P.S.Prudnikov for comments on the paper.



\begin{thebibliography}{99}

\bibitem{remiz:1}
K.-J. Engel, R. Nagel. One-Parameter Semigroups for Linear Evolution Equations. --- Springer, 2000.

\bibitem{remiz:2}
V.I.Bogachev, O.G.Smolyanov. Real and functional analysis (University course), 2nd edition. --- Regular and chaotic dynamics, Izhevsk, 2011

\bibitem{Chernoff}
Paul R. Chernoff, Note on product formulas for operator semigroups. // J. Funct. Anal. 2:2 (1968) 238-242.

\bibitem{remiz:4}
I.D. Remizov. Quasi-Feynman formulas --- a method of obtaining the evolution operator for the Schr\"odinger equation. // Journal of Functional Analysis 270:12 (2016) 4540-4557


\bibitem{OSS2012} 
Yu. N. Orlov, V. Zh. Sakbaev, O. G. Smolyanov. Rate of convergence of Feynman approximations of semigroups generated by the oscillator Hamiltonian// TMF, 172:1 (2012),  122-137 

\bibitem{Remizov-NNGU2018}
 I.\,D.~Remizov. On estimation of error in approximations provided by Chernoff's product formula// International Conference "ShilnikovWorkshop-2018", Lobachevsky State University of Nizhny Novgorod (Russia), book of abstracts, pp.38-41 (2018)


\bibitem{Rem-add-2}
I.D. Remizov. Approximations to the solution of Cauchy problem for a linear evolution equation via the space shift operator (second-order equation example). // Applied Mathematics and Computaton 328 (2018), 243-246.

\bibitem{remiz:3}
I.D. Remizov. Feynman and Quasi-Feynman Formulas for Evolution Equations. // Doklady Mathematics, 96:2 (2017), 433-437

\bibitem{R-Starod}
I. D. Remizov, M. F. Starodubtseva. Quasi-Feynman Formulas providing Solutions of Multidimensional Schr\"odinger Equations with Unbounded Potential.// Math. Notes, 104:5 (2018), 767-772

\bibitem{Loboda}
A. A. Loboda. The Doss Method for the Stochastic Schr\"odinger-Belavkin Equation// Math. Notes, 106:2 (2019), 303-307

\bibitem{Zag}
V. Zagrebnov. Notes on the Chernoff Product Formula// Cornell University library archive, arXiv:1911.09480 (2019)

\bibitem{Gom}
A. Gomilko, S. Kosowicz, Yu. Tomilov. A general approach to approximation theory of operator semigroups//Cornell University library archive, arXiv:1801.06749 (2018)

\end{thebibliography}
\end{document}